\documentclass[a4paper,12pt,reqno]{amsart}
\usepackage{amssymb,amsthm}
\usepackage{ifthen}
\usepackage{graphicx}
\usepackage{float}

\theoremstyle{plain}

\newtheorem{thm}{Theorem}[section]
\newtheorem*{thm*}{Theorem}

\numberwithin{equation}{section}

\newtheorem{cor}{Corollary}[section]

\newtheorem{prop}{Proposition}[section]
\newtheorem{rem}{Remark}[section]

\newtheorem{defn}{Definition}[section]

\theoremstyle{definition}
\newcounter {own}
\def\theown {\thesection  .\arabic{own}}

\newenvironment{pf}[1][]{%
 \vskip 3mm
 \noindent
 \ifthenelse{\equal{#1}{}}%
  {{\slshape Proof. }}%
  {{\slshape #1.} }%
 }%
{\qed\bigskip}

\newcounter{alphabet}
\newcounter{tmp}

\newcommand{\ds}{\displaystyle}

\newcounter{minutes}
\divide\time by 60
\newcounter{hours}
\multiply\time by 60 \addtocounter{minutes}{-\time}

\begin{document}
\bibliographystyle{amsplain}
\title{Separation of zeros and a Hermite interpolation based frame algorithm for band limited functions}

\thanks{
File:~\jobname .tex,
          printed: 2013-07-30,
          \thehours.\ifnum\theminutes<10{0}\fi\theminutes}

\author{A. Antony Selvan}

\address{A. Antony Selvan, Department of Mathematics,
Indian Institute of Technology Madras, Chennai--600 036, India.}
\email{antonyaans@gmail.com}
\author{R. Radha $^\dagger$}
\address{R. Radha, Department of Mathematics,
Indian Institute of Technology Madras, Chennai--600 036, India.}
\email{radharam@iitm.ac.in}
\keywords{band limited functions, Bernstein's inequality, frames, Hermite interpolation, nonuniform sampling, Wirtinger-Sobolev inequality.\\
$^\dagger$ {\tt Corresponding author}
}
\maketitle
\pagestyle{myheadings}
\markboth{A. Antony Selvan and R. Radha}{Separation of zeros and a Hermite interpolation based frame algorithm for band limited functions}

\begin{abstract}
It is shown that if a non-zero function $f\in B_\sigma$  has infinitely many double zeros on the real axis, then there exists at least one pair of consecutive zeros whose distance apart is greater than $\dfrac{\pi}{\sigma}\tau^{1/4}$, $\tau\approx5.0625$. A frame algorithm  is provided for reconstructing a function $f\in B_\sigma$ from its nonuniform samples $\{f^{(j)}(x_i):j=0,1,\dots, k-1, i\in\mathbb{Z}\}$ with maximum gap condition, $\sup\limits_i(x_{i+1}-x_i)=\delta<\dfrac{1}{\sigma}c_k^{1/2k}$, where $c_k$ is a Wirtinger-Sobolev constant, using Hermite interpolation.
\end{abstract}
\section{Introduction}
Let $B_\sigma$ denote the space of all $\sigma$-band limited functions, \textit{i.e.,}
\begin{eqnarray*}
B_\sigma=\left\{f\in L^2(\mathbb{R}): \textrm{supp}\widehat{f}\subseteq [-\sigma,\sigma]\right\}.
\end{eqnarray*}
Here $\widehat{f}$ denotes the Fourier transform of $f$, defined by $\widehat{f}(\xi)=\int_{-\infty}^\infty f(x)e^{-ix\xi}\mathrm{d}x$. Then it follows from the well known theorem of Paley-Wiener that a function $f\in B_\sigma$ if and only if  $f$ can be extended as an entire function  of exponential type $\leq\sigma$.
It is well known that $B_\sigma$ is a reproducing kernel Hilbert space with reproducing kernel $K(x,y)=\dfrac{\sin \sigma\left(x-y\right)}{\sigma \left(x-y\right)}$. The classical Shannon's sampling theorem states
that every $f\in B_\sigma$ can be reconstructed from the sampling formula
\begin{eqnarray*}
f(x)=\ds\sum_{k\in\mathbb{Z}}f\left(\dfrac{k\pi}{\sigma}\right)\dfrac{\sin \sigma\left(x-k\pi/\sigma\right)}{\sigma \left(x-k\pi/\sigma\right)}.
\end{eqnarray*}
In \cite{Walker} Walker proved the following result. ``Let $f$ be an entire function of exponential type $\leq \sigma$ which is square integrable on the real axis. If $f$ has infinitely many zeros on the real axis, then there exists at least one pair of consecutive zeros whose distance apart is at least $\dfrac{\pi}{\sigma}$''. In this paper we show that if a non-zero function $f\in B_\sigma$  has infinitely many double zeros on the real axis, then there exists at least one pair of consecutive zeros whose distance apart is greater than $\dfrac{\pi}{\sigma}\tau^{1/4}$, $\tau\approx5.0625$.

The numerical aspects of nonuniform sampling was studied by Grochenig in \cite{Grochenig}. In particular, in that paper, he gave an iterative reconstruction algorithm for a band limited function from its nonuniform samples and discussed its stability as well as its rate of convergence. In \cite{Raza}, Razafinjatovo obtained a frame algorithm for reconstructing a function $f\in B_\sigma$ from its nonuniform samples $\{f^{(j)}(x_i):j=0,1,\dots, k-1, i\in\mathbb{Z}\}$ with maximum gap condition, namely  $\sup\limits_i(x_{i+1}-x_i)=\delta<\dfrac{2}{\sigma}((k-1)!\sqrt{(2k-1)2k})^{1/k}$ using Taylor's polynomial approximation. He has also provided a result which gives an improvement for the maximum gap condition for $k=2$.

In this paper, we provide a frame algorithm for reconstructing a function $f\in B_\sigma$ from its nonuniform samples $\{f^{(j)}(x_i):j=0,1,\dots, k-1, i\in\mathbb{Z}\}$ with maximum gap condition, $\delta<\dfrac{1}{\sigma}c_k^{1/2k}$, where $c_k$ is a Wirtinger-Sobolev constant, using Hermite interpolation. We wish to remark that the maximum gap condition of the current paper is an improvement of the maximum gap condition given in \cite{Raza}. In fact, later, we provide a numerical comparison between the maximum gap conditions of \cite{Raza} and the current paper for various values of $k$. Further, our idea of using Hermite interpolation polynomial instead of the already existing method using Taylor series expansion leads to better rate of convergence, which is mentioned in Remark 2.1.

It should also be noted that if the sample points $\{x_i\}$ satisfy  $\sup\limits_i(x_{i+1}-x_i)=\delta<\dfrac{\pi}{\sigma}\tau^{1/4}, \tau\approx5.0625$, then from Theorem \ref{pap3thm2.1} and  \ref{pap3thm2.3}, it follows that one can reconstruct functions  $f\in B_\sigma$ uniquely from $\{f(x_i), f'(x_i)\}$. In addition, as a consequence of  Theorem \ref{pap3thm2.1}, one can get the uniqueness result even if the equality holds.

In order to prove our main results, we make use of the following terminology and some inequalities.
\begin{defn}
A sequence of vectors $\{f_n:~n\in\mathbb{Z}\}$ in  a separable
Hilbert space $\mathcal{H}$ is said to be a {\it frame} if there
exist constants $0< m \leq M<\infty$ such that
\begin{equation}
m\|f\|_\mathcal{H}^2\leq\displaystyle\sum_{n\in\mathbb{Z}}|\langle
f,f_n\rangle_\mathcal{H}|^2\leq M\|f\|_\mathcal{H}^2,
\end{equation}
for every $f\in\mathcal{H}$.
\end{defn}

\begin{thm}[Bernstein's inequality]
If $f\in B_\sigma$, then
\begin{eqnarray}
\|f^{(k)}\|_2\leq \sigma^k\|f\|_2.
\end{eqnarray}
\end{thm}
\begin{thm}[Arthurs-Anderson-Hall inequality]\label{ArAnHaineq}$(\mathop{\rm cf.~}$\cite{Horgan}$)$.
Let $f$ be a complex valued  function defined on the interval $\left[a,b\right]$. If $f\in C^2\left[a,b\right]$ with $f(a)=f(b)=f'(a)=f'(b)=0$, then \begin{eqnarray}
\int\limits_a^b|f(x)|^2~\mathrm{d}x\leq \dfrac{1}{\tau}\Big(\dfrac{b-a}{\pi}\Big)^4\int\limits_a^b|f^{\prime\prime}(x)|^2~\mathrm{d}x,
\end{eqnarray}
where $\tau$ is the smallest root of the equation
\begin{eqnarray*}
\tanh\Big(\dfrac{\pi \tau^{1/4}}{2}\Big)+\tan\Big(\dfrac{\pi \tau^{1/4}}{2}\Big)=0.
\end{eqnarray*}
The equality holds iff
\begin{eqnarray*}
f(x)=c\left\{\cosh\mu\dfrac{x-a}{b-a}-\cos\mu\dfrac{x-a}{b-a}-\dfrac{\cosh\mu-cos\mu}{\sinh\mu-sin\mu}
\left(\sinh\mu\dfrac{x-a}{b-a}-\sin\mu\dfrac{x-a}{b-a}\right)\right\},
\end{eqnarray*}
$c\in\mathbb{C}$ with $\mu=\pi \tau^{1/4}$. The value of $\tau$ is approximately $5.0625$.
\end{thm}
\begin{thm}[Wirtinger-Sobolev inequality]\label{pap2WirSob}$(\mathop{\rm cf.~}$\cite{BoWi}$)$.
Let $f$ be a complex valued  function defined on the interval $\left[a,b\right]$.
If $f\in C^r\left[a,b\right]$ with $f^{(l)}(a)=f^{(l)}(b)=0$, $0\leq l\leq r-1$, then
\begin{eqnarray}\label{pap3eqn1.4}
\int\limits_a^b|f(x)|^2~\mathrm{d}x\leq \dfrac{(b-a)^{2r}}{c_r}\int\limits_a^b|f^{(r)}(x)|^2~\mathrm{d}x,
\end{eqnarray}
where $c_r$ is the minimal eigenvalue of the boundary value problem
\begin{eqnarray*}
(-1)^r u^{(2r)}(x)=\lambda u(x),~~ x\in[0,1],\\
u^{(k)}(0)=u^{(k)}(1)=0, ~~0\leq k\leq r-1,
\end{eqnarray*}
$u\in C^{2r}[0,1]$.
\end{thm}

\begin{prop}$(\mathop{\rm cf.~}$\cite{FeiGr},\cite{Raza}$)$.\label{pap2prop2.1}
Let $A$ be a bounded operator on a Hilbert space $\mathcal{H}$ that satisfies
\begin{align*}
\|f-Af\|_\mathcal{H}\leq C\|f\|_\mathcal{H},
\end{align*}
for every $f\in \mathcal{H}$ and for some $C$, $0<C<1$.
Then $A$ is invertible on $\mathcal{H}$ and $f$ can be recovered from $Af$ by the following iteration algorithm.
Setting
\begin{eqnarray*}
f_0&=& A f ~and \\
f_{n+1}&=&f_n+A(f-f_n), ~n\geq 0,
\end{eqnarray*}
we have $\lim\limits_{n\to\infty}f_n=f$. The error estimate after $n$ iterations is
\begin{align*}
\|f-f_n\|_\mathcal{H}\leq C^{n+1}\|f\|_\mathcal{H}.
\end{align*}
\end{prop}
\begin{thm}[Hermite Interpolation Formula] $(\mathop{\rm cf.~}$\cite{Spitzbart}$)$.
Let $f\in C^r\left[a,b\right]$ and $\xi,\eta\in[a,b]$. Then the Hermite interpolation polynomial $H_{2r+1}(x)$ of degree $2r+1$ such that
$H_{2r+1}^{(j)}(y)=f^{(j)}(y)$, for $y=\xi,\eta$, $0\leq j\leq r$, is given by
\begin{eqnarray}
H_{2r+1}(\xi,\eta,f;x)=\ds\sum\limits_{k=0}^{r}A_{0k}(x)f^{(k)}(\xi)+\ds\sum\limits_{k=0}^{r}A_{1k}(x)f^{(k)}(\eta),
\end{eqnarray}
where
\begin{eqnarray*}
A_{0k}(x)&=&(x-\eta)^{r+1}\dfrac{(x-\xi)^k}{k!}\ds\sum\limits_{s=0}^{r-k}\dfrac{1}{s!}g_0^{(s)}(\xi)(x-\xi)^s,\\
A_{1k}(x)&=&(x-\xi)^{r+1}\dfrac{(x-\eta)^k}{k!}\ds\sum\limits_{s=0}^{r-k}\dfrac{1}{s!}g_1^{(s)}(\eta)(x-\eta)^s,\\
g_0(x)&=&(x-\eta)^{-(r+1)},\\
g_1(x)&=&(x-\xi)^{-(r+1)}.
\end{eqnarray*}
\end{thm}
\section{The Main Results}

First we observe that Wirtinger-Sobolev inequality (Theorem \eqref{pap2WirSob}) is still true if the right hand side of \eqref{pap3eqn1.4} is replaced by
\begin{eqnarray*}
\dfrac{1}{c_r^2}(b-a)^{4r}\int\limits_a^b|f^{(2r)}(x)|^2~\mathrm{d}x,
\end{eqnarray*}
$f\in C^{2r}\left[a,b\right]$ with $f^{(l)}(a)=f^{(l)}(b)=0$, $0\leq l\leq r-1$. In other words,
\begin{eqnarray}\label{pap3eqn2.1}
\|f\|_{L^2[a,b]}^2\leq \dfrac{1}{c_r^2}(b-a)^{4r}\|f^{(2r)}\|_{L^2[a,b]}^2.
\end{eqnarray}
In fact, we shall establish \eqref{pap3eqn2.1} by assuming $f$ is real valued. The complex case will follow by taking $f=u+iv$. We know that $A^2-B^2=(A+B)^2-2(AB+B^2)\geq-2(AB+B^2)$ for real numbers $A$ and $B$. Therefore,
\begin{eqnarray*}
\int\limits_a^b \left[\dfrac{(b-a)^{4r}}{c_r^2}f^{(2r)}(x)^2-f^2(x)\right]\mathrm{d}x
&\geq& -2\int\limits_a^b \left[\dfrac{(b-a)^{2r}}{c_r}f^{(2r)}(x)f(x)+f^2(x)\right]\mathrm{d}x.
\end{eqnarray*}
Using Bernoulli's formula (repeated integration by parts) and the fact that $f^{(l)}(a)=f^{(l)}(b)=0$, for $0\leq l\leq r-1$, $r$ odd, one can easily show that
\begin{eqnarray*}
-2\int\limits_a^b \left[\dfrac{(b-a)^{2r}}{c_r}f^{(2r)}(x)f(x)+f^2(x)\right]\mathrm{d}x
&\geq&2\int\limits_a^b \Big[\dfrac{(b-a)^{2r}}{c_r}f^{(r)}(x)^2\Big]\mathrm{d}x-2 \int\limits_a^b f^2(x)\mathrm{d}x,\\
&\geq&0,
\end{eqnarray*}
by \eqref{pap3eqn1.4}, thus proving \eqref{pap3eqn2.1}, for $r$ odd.

We know that $A^2-B^2=(A-B)^2+2(AB-B^2)\geq2(AB-B^2)$ for real numbers $A$ and $B$. Therefore,
\begin{eqnarray*}
\int\limits_a^b \left[\dfrac{(b-a)^{4r}}{c_r^2}f^{(2r)}(x)^2-f^2(x)\right]\mathrm{d}x
&\geq& 2\int\limits_a^b \left[\dfrac{(b-a)^{2r}}{c_r}f^{(2r)}(x)f(x)-f^2(x)\right]\mathrm{d}x.
\end{eqnarray*}
Then proceeding as before we obtain \eqref{pap3eqn2.1}, for $r$ even.

\begin{thm}\label{pap3thm2.1}
If a non-zero function $f\in B_\sigma$  has infinitely many double zeros on the real axis, then there exists at least one pair of consecutive zeros whose distance apart is greater than $\dfrac{\pi}{\sigma}\tau^{1/4}$, where $\tau\approx5.0625$.
\end{thm}
\begin{pf}
Let a non zero function $f\in B_\sigma$ have infinitely many double zeros $x_j$'s on the real line such that $x_j<x_{j+1}$, $j\in\mathbb{Z}$ and $\bigcup\limits_{j\in\mathbb{Z}}[x_j,x_{j+1}]=\mathbb{R}$.
If possible, there exists $M\leq\dfrac{\pi}{\sigma}\tau^{1/4}$ such that $x_{j+1}-x_j\leq M$, for every $j$. Since $f(x_j)=f'(x_j)=f(x_{j+1})=f'(x_{j+1})=0$, for every $j$, Arthurs-Anderson-Hall inequality yields
\begin{eqnarray}\label{pap3eqn2.3}
\ds\int\limits_{x_j}^{x_{j+1}}|f(x)|^2~\mathrm{d}x<\dfrac{1}{\tau}\left(\dfrac{x_{j+1}-x_j}{\pi}\right)^4\ds\int\limits_{x_j}^{x_{j+1}}|f^{\prime\prime}(x)|^2~\mathrm{d}x.
\end{eqnarray}
Notice that the inequality is strict; Otherwise if the equality holds, then, by uniqueness theorem,
\begin{eqnarray*}
f(x)=c\left\{\cosh\mu\dfrac{x-a}{b-a}-\cos\mu\dfrac{x-a}{b-a}-\dfrac{\cosh\mu-cos\mu}{\sinh\mu-sin\mu}
\left(\sinh\mu\dfrac{x-a}{b-a}-\sin\mu\dfrac{x-a}{b-a}\right)\right\},
\end{eqnarray*}
$c\in\mathbb{C}-\{0\}$ with $a=x_j$ , $b=x_{j+1}$, $\mu=\pi\tau^{1/4}$, $x\in\mathbb{R}$. Clearly, $f\notin L^2(\mathbb{R})$, which is impossible.
Summing over all $j$ in \eqref{pap3eqn2.3}, we get
\begin{eqnarray*}
\ds\int\limits_{\mathbb{R}}|f(x)|^2~\mathrm{d}x
&<&\ds\sum\limits_{j}\dfrac{1}{\tau}\left(\dfrac{x_{j+1}-x_j}{\pi}\right)^4\ds\int\limits_{x_j}^{x_{j+1}}|f^{\prime\prime}(x)|^2~\mathrm{d}x\\
&\leq&\dfrac{1}{\tau}\dfrac{M^4}{\pi^4}\ds\int\limits_{\mathbb{R}}|f^{\prime\prime}(x)|^2~\mathrm{d}x.
\end{eqnarray*}
Taking square root on both sides, we get
\begin{eqnarray}\label{pap3eqn2.4}
\|f\|_2<\dfrac{1}{\sqrt{\tau}}\dfrac{M^2}{\pi^2}\|f^{\prime\prime}\|_2.
\end{eqnarray}
On the other hand, by Bernstein's inequality,
\begin{eqnarray}\label{pap3eqn2.5}
\|f^{\prime\prime}\|_2\leq\sigma^2~\|f\|_2.
\end{eqnarray}
Combining $\eqref{pap3eqn2.4}$ and $\eqref{pap3eqn2.5}$, we get $M>\dfrac{\pi}{\sigma}\tau^{1/4}$ which is a contradiction.
\end{pf}
\begin{cor}
Let $f\in B_\sigma$ satisfy $f(x_j)=f'(x_j)=0$, for all $j\in\mathbb{Z}$. If $\sup\limits_{j}(x_{j+1}-x_j)\leq\dfrac{\pi}{\sigma}\tau^{1/4}$, then $f\equiv 0$.
\end{cor}
Consider the operator $P:L^2(\mathbb{R})\to B_\sigma$ by
\begin{eqnarray}\label{pap2eqn2.6}
(Pf)(x):=\langle f, K_x\rangle,
\end{eqnarray}
where $K_x(t)=\dfrac{\sin \sigma\left(t-x\right)}{\sigma \left(t-x\right)}$. Then $P$ is an orthogonal projection of $L^2(\mathbb{R})$ onto $B_\sigma$.
Now assume that $f$ and its first $k-1$ derivatives $f', \dots,f^{(k-1)}$ are sampled at a sequence $(x_i)_{i\in\mathbb{Z}}$.
Define the approximation operator for $f\in B_\sigma$
\begin{align*}
Af=P\left(\ds\sum\limits_{i\in\mathbb{Z}}H_{2k-1}(x_i,x_{i+1},f;\cdot)\chi_{[x_i,x_{i+1}]}\right).
\end{align*}
Since $H_{2r+1}(\xi,\eta,\alpha f+g;x)=\alpha H_{2r+1}(\xi,\eta,f;x)+H_{2r+1}(\xi,\eta,g;x)$ for $\alpha\in\mathbb{C}$, the operator $A$ is linear. Since  $f=Pf=P\left(\ds\sum\limits_{i\in\mathbb{Z}}f\chi_{[x_i,x_{i+1}]}\right)$, and the characteristic functions $\chi_{[x_i,x_{i+1}]}$ have mutually disjoint support it can be easily shown that
\begin{eqnarray}
\|f-Af\|_2^2&\leq&\ds\sum\limits_{i\in\mathbb{Z}}\int\limits_{x_i}^{x_{i+1}}|f(x)-H_{2k-1}(x_i,x_{i+1},f;x)|^2~\mathrm{d}x,
\end{eqnarray}
where $H_{2k-1}(x_i,x_{i+1},f;\cdot)$ denotes the Hermite interpolation of $f$ in the interval $[x_i,x_{i+1}]$.
Let $\sup\limits_i(x_{i+1}-x_i)=\delta<\dfrac{1}{\sigma}c_{k}^{\frac{1}{2k}}$. Using \eqref{pap3eqn2.1}, we get
\begin{eqnarray}
\|f-Af\|_2^2
&\leq&\ds\sum\limits_{i\in\mathbb{Z}}\dfrac{\delta^{4k}}{c_{k}^2}
\int\limits_{x_i}^{x_{i+1}}|f^{(2k)}(x)|^2~\mathrm{d}x\nonumber\\
&=&\dfrac{\delta^{4k}}{c_{k}^2}\|f^{(2k)}\|_2^2\nonumber\\
&\leq&\dfrac{\delta^{4k}\sigma^{4k}}{c_{k}^2}\|f\|_2^2,
\end{eqnarray}
using Bernstein's inequality. As $\|Af\|_2\leq \|f-Af\|_2+\|f\|_2$, it follows from the inequality \eqref{pap3eqn2.1} that the operator $A$ is a bounded.  Now, if $\delta<\dfrac{1}{\sigma}c_{k}^{\frac{1}{2k}}$, then $\dfrac{\delta^{4k}\sigma^{4k}}{c_{k}^2}<1$. Thus  we can obtain the following result as a corollary of Proposition \ref{pap2prop2.1}.
\begin{thm}\label{pap3thm2.2}
Assume that $f$ and its first $k-1$ derivatives $f', \dots,f^{(k-1)}$ are sampled at a sequence $(x_i)_{i\in\mathbb{Z}}$. If $\delta<\dfrac{1}{\sigma}c_{k}^{\frac{1}{2k}}$, then  any $f\in B_\sigma$ can be reconstructed from the sample values $\{f^{(j)}(x_i):j=0,1,\dots,k-1,i\in\mathbb{Z}\}$ using the following iteration algorithm.
Set
\begin{eqnarray*}
f_0&=&A f=P\left(\ds\sum\limits_{i\in\mathbb{Z}}H_{2k-1}(x_i,x_{i+1},f;\cdot)\chi_{[x_i,x_{i+1}]}\right),\\
f_{n+1}&=&f_n+A(f-f_n),~ n\geq 0,
\end{eqnarray*}
where $H_{2k-1}(x_i,x_{i+1},f;\cdot)$ denotes the Hermite interpolation of $f$ in the interval $[x_i,x_{i+1}]$.
Then we have $\lim\limits_{n\to\infty}f_n=f$. The error estimate after $n$ iterations becomes
\begin{eqnarray*}
\|f-f_n\|_2&\leq&\left(\dfrac{(\delta\sigma)^{2k}}{c_k}\right)^{(n+1)}\|f\|_2.
\end{eqnarray*}
\end{thm}
\begin{rem}
The above inequality shows that the rate of convergence of the current paper is better than that of \cite{Raza}.
\end{rem}

In order to present frame algorithm for reconstructing a function $f\in B_\sigma$ from  its nonuniform samples we need the following notations.\\
Let $c_{i,l}=\ds\int\limits_{x_i}^{x_{i+1}}\dfrac{(x-x_{i+1})^{2l}}{l!^2}~\mathrm{d}x$. This can also be written as
\begin{eqnarray*}
c_{i,l}=\dfrac{(x_{i+1}-x_i)^{2l+1}}{(2l+1)l!^2}=\ds\int\limits_{x_i}^{x_{i+1}}\dfrac{(x-x_{i})^{2l}}{l!^2}~\mathrm{d}x.
\end{eqnarray*}
Let $g_0(x)=(x-\eta)^{-k}$ and $g_1(x)=(x-\xi)^{-k}$. Then
\begin{eqnarray*}
g_0^{(s)}(x)&=&(-1)^s k (k+1)\dots(k-1+s)(x-\eta)^{-(k+s)}~\textrm{and}\\
g_1^{(s)}(x)&=&(-1)^s k (k+1)\dots(k-1+s)(x-\xi)^{-(k+s)}.
\end{eqnarray*}
\begin{thm}\label{pap3thm2.3}
If $\delta<\dfrac{1}{\sigma}c_{k}^{\frac{1}{2k}}$, then for every $f\in B_\sigma$, we have
\begin{eqnarray}\label{pap3eqn2.7}
A\| f\|_2^2&\leq&\ds\sum\limits_{i\in\mathbb{Z}}
\ds\sum\limits_{l=0}^{k-1}|f^{(l)}(x_i)|^2(c_{i,l}+c_{i-1,l})\leq B\| f\|_2^2,
\end{eqnarray}
where
$A=\left(1-\dfrac{(\delta\sigma)^{2k}}{c_k}\right)^{2}\dfrac{1}{2k C(k)}$,
$B=2\left(\ds\sum\limits_{l=0}^{k-1}\dfrac{(\delta\sigma)^{2l}}{l!^2}\right)e^{\delta^2+\sigma^2}$, and\\
$C(k)=\left[\ds\sum\limits_{s=0}^{k-1}{k+s-1\choose s} \right]^2$.
\end{thm}
\begin{pf}
Recall $Af=P\left(\ds\sum\limits_{i\in\mathbb{Z}}H_{2k-1}(x_i,x_{i+1},f;\cdot)\chi_{[x_i,x_{i+1}]}\right)$. Then
\begin{eqnarray}\label{pap3eqn2.8}
\|f\|_2^2&=&\|A^{-1}Af\|_2^2\leq\|A^{-1}\|^2\|Af\|_2^2\nonumber\\
&\leq&(1-\|I-A\|)^{-2}\|Af\|_2^2\nonumber\\
&\leq&\left(1-\dfrac{(\delta\sigma)^{2k}}{c_k}\right)^{-2}\|Af\|_2^2.
\end{eqnarray}
We now estimate the value of $\|Af\|_2$.
\begin{eqnarray}
\|A f\|_2^2
&\leq&\left\|\ds\sum\limits_{i\in\mathbb{Z}}H_{2k-1}(x_i,x_{i+1},f;\cdot)\chi_{[x_i,x_{i+1}]}\right\|_2^2\nonumber\\
&=&\ds\int\limits_{\mathbb{R}}\left|\ds\sum\limits_{i\in\mathbb{Z}}H_{2k-1}(x_i,x_{i+1},f;x)\chi_{[x_i,x_{i+1}]}(x)\right|^2~\mathrm{d}x\nonumber\\
&\leq&\ds\sum\limits_{i\in\mathbb{Z}}\int\limits_{x_i}^{x_{i+1}}|H_{2k-1}(x_i,x_{i+1},f;x)|^2~\mathrm{d}x\nonumber\\
&=&\ds\sum\limits_{i\in\mathbb{Z}}\int\limits_{x_i}^{x_{i+1}}
\left|\ds\sum\limits_{l=0}^{k-1}A_{0l}(x)f^{(l)}(x_i)+\ds\sum\limits_{l=0}^{k-1}A_{1l}(x)f^{(l)}(x_{i+1})\right|^2~\mathrm{d}x\nonumber\\
&\leq&2\left\{\ds\sum\limits_{i\in\mathbb{Z}}\int\limits_{x_i}^{x_{i+1}}
\left|\ds\sum\limits_{l=0}^{k-1}A_{0l}(x)f^{(l)}(x_i)\right|^2+\left|\ds\sum\limits_{l=0}^{k-1}A_{1l}(x)f^{(l)}(x_{i+1})\right|^2~\mathrm{d}x\right\}\nonumber\\
&\leq&2k\left\{\ds\sum\limits_{i\in\mathbb{Z}}\int\limits_{x_i}^{x_{i+1}}
\ds\sum\limits_{l=0}^{k-1}|A_{0l}(x)|^2|f^{(l)}(x_i)|^2+\ds\sum\limits_{l=0}^{k-1}|A_{1l}(x)|^2|f^{(l)}(x_{i+1})|^2~\mathrm{d}x\right\}\nonumber\\
&=&2k\left\{\ds\sum\limits_{i\in\mathbb{Z}}
\ds\sum\limits_{l=0}^{k-1}|f^{(l)}(x_i)|^2\int\limits_{x_i}^{x_{i+1}}|A_{0l}(x)|^2\mathrm{d}x
+\ds\sum\limits_{i\in\mathbb{Z}}\ds\sum\limits_{l=0}^{k-1}|f^{(l)}(x_{i+1})|^2\int\limits_{x_i}^{x_{i+1}}|A_{1l}(x)|^2~\mathrm{d}x\right\}.\nonumber\\
\end{eqnarray}
Let $I_1=\ds\int\limits_{x_i}^{x_{i+1}}|A_{0l}(x)|^2~\mathrm{d}x$ and $I_2=\ds\int\limits_{x_i}^{x_{i+1}}|A_{1l}(x)|^2~\mathrm{d}x$. Now,
\begin{eqnarray}
I_1&=&\int\limits_{x_i}^{x_{i+1}}
\left[\ds\sum\limits_{s=0}^{k-1-l}\dfrac{1}{s!}g_0^{(s)}(x_i)(x-x_{i})^s\right]^2(x-x_{i+1})^{2k}\dfrac{(x-x_{i})^{2l}}{l!^2}~\mathrm{d}x\nonumber\\
&\leq&\max\limits_{[x_i,x_{i+1}]}\left[\ds\sum\limits_{s=0}^{k-1-l}\dfrac{1}{s!}g_0^{(s)}(x_i)(x-x_{i})^s\right]^2
\max\limits_{[x_i,x_{i+1}]}(x-x_{i+1})^{2k}\times c_{i,l}\nonumber\\
&=&(x_i-x_{i+1})^{2k}\times c_{i,l}\nonumber\\
&&\times\max\limits_{[x_i,x_{i+1}]}\left[\ds\sum\limits_{s=0}^{k-1-l}\dfrac{(-1)^s}{s!} k (k+1)\dots(k-1+s)(x_i-x_{i+1})^{-(k+s)}(x-x_{i})^s\right]^2\nonumber
\end{eqnarray}
\begin{eqnarray}
&\leq&c_{i,l}\max\limits_{[x_i,x_{i+1}]}\left[\ds\sum\limits_{s=0}^{k-1-l}\dfrac{1}{s!} k (k+1)\dots(k-1+s)(x_{i+1}-x_{i})^{-s}(x-x_{i})^s\right]^2\nonumber\\
&=& c_{i,l}\left[\ds\sum\limits_{s=0}^{k-1-l}\dfrac{1}{s!} k (k+1)\dots(k-1+s)\right]^2\nonumber\\
&\leq& c_{i,l}\left[\ds\sum\limits_{s=0}^{k-1}{k+s-1\choose s}\right]^2.\nonumber
\end{eqnarray}
Thus $I_1\leq C(k)c_{i,l}$. Similarly, $I_2\leq C(k)c_{i,l}$. Hence
\begin{eqnarray}
\|A f\|_2^2&\leq&2k C(k)\left\{\ds\sum\limits_{i\in\mathbb{Z}}
\ds\sum\limits_{l=0}^{k-1}|f^{(l)}(x_i)|^2c_{i,l}
+\ds\sum\limits_{i\in\mathbb{Z}}\ds\sum\limits_{l=0}^{k-1}|f^{(l)}(x_{i+1})|^2c_{i,l}\right\}\nonumber\\
&\leq&2k C(k)\ds\sum\limits_{i\in\mathbb{Z}}
\ds\sum\limits_{l=0}^{k-1}|f^{(l)}(x_i)|^2(c_{i,l}+c_{i-1,l}).
\end{eqnarray}
Therefore, \eqref{pap3eqn2.8} becomes
\begin{eqnarray}
\| f\|_2^2&\leq&\left(1-\dfrac{(\delta\sigma)^{2k}}{c_k}\right)^{-2}2k C(k)\ds\sum\limits_{i\in\mathbb{Z}}
\ds\sum\limits_{l=0}^{k-1}|f^{(l)}(x_i)|^2(c_{i,l}+c_{i-1,l}).
\end{eqnarray}
Thus, we obtain the LHS of \eqref{pap3eqn2.7}. For any $f\in B_\sigma$, we have
\begin{eqnarray}
f^{(l)}(x_i)=\ds\sum\limits_{\nu=0}^{\infty}f^{(l+\nu)}(x)\dfrac{(x_i-x)^{\nu}}{\nu!}.
\end{eqnarray}
Hence,
\begin{eqnarray}\label{pap3eqn2.13}
\ds\int\limits_{x_i}^{x_{i+1}}|f^{(l)}(x_i)|^2~\mathrm{d}x
&=&\ds\int\limits_{x_i}^{x_{i+1}}\left|\ds\sum\limits_{\nu=0}^{\infty}f^{(l+\nu)}(x)\dfrac{(x_i-x)^{\nu}}{\nu!}\right|^2~\mathrm{d}x\nonumber\\
&\leq&\ds\int\limits_{x_i}^{x_{i+1}}
\left(\ds\sum\limits_{\nu=0}^{\infty}\dfrac{\left|f^{(l+\nu)}(x)\right|^2}{\nu!}\right)
\left(\ds\sum\limits_{\nu=0}^{\infty}\dfrac{(x_i-x)^{2\nu}}{\nu!}\right)~\mathrm{d}x\nonumber\\
&\leq& e^{\delta^2}\ds\sum\limits_{\nu=0}^{\infty}\ds\int\limits_{x_i}^{x_{i+1}}
\dfrac{\left|f^{(l+\nu)}(x)\right|^2}{\nu!}~\mathrm{d}x.
\end{eqnarray}
Now,
\begin{eqnarray}
\ds\sum\limits_{i\in\mathbb{Z}}\ds\sum\limits_{l=0}^{k-1}|f^{(l)}(x_i)|^2c_{i,l}
&=& \ds\sum\limits_{i\in\mathbb{Z}}\ds\sum\limits_{l=0}^{k-1}|f^{(l)}(x_i)|^2\ds\int\limits_{x_i}^{x_{i+1}}\dfrac{(x-x_{i+1})^{2l}}{l!^2}~\mathrm{d}x\nonumber\\
&=& \ds\sum\limits_{l=0}^{k-1}\ds\sum\limits_{i\in\mathbb{Z}}\ds\int\limits_{x_i}^{x_{i+1}}|f^{(l)}(x_i)|^2\dfrac{(x-x_{i+1})^{2l}}{l!^2}~\mathrm{d}x\nonumber\\
&\leq& \ds\sum\limits_{l=0}^{k-1}\ds\sum\limits_{i\in\mathbb{Z}}\dfrac{\delta^{2l}}{l!^2}\ds\int\limits_{x_i}^{x_{i+1}}|f^{(l)}(x_i)|^2~\mathrm{d}x\nonumber\\
&\leq& \ds\sum\limits_{l=0}^{k-1}\dfrac{\delta^{2l}}{l!^2}\ds\sum\limits_{i\in\mathbb{Z}}
e^{\delta^2}\ds\sum\limits_{\nu=0}^{\infty}\ds\int\limits_{x_i}^{x_{i+1}}
\dfrac{\left|f^{(l+\nu)}(x)\right|^2}{\nu!}~\mathrm{d}x,\nonumber
\end{eqnarray}
using \eqref{pap3eqn2.13}. Hence
\begin{eqnarray}
\ds\sum\limits_{i\in\mathbb{Z}}\ds\sum\limits_{l=0}^{k-1}|f^{(l)}(x_i)|^2c_{i,l}
&\leq& \ds\sum\limits_{l=0}^{k-1}\dfrac{\delta^{2l}}{l!^2}e^{\delta^2}
\ds\sum\limits_{\nu=0}^{\infty}\ds\sum\limits_{i\in\mathbb{Z}}\ds\int\limits_{x_i}^{x_{i+1}}
\dfrac{\left|f^{(l+\nu)}(x)\right|^2}{\nu!}~\mathrm{d}x\nonumber\\
&=& \ds\sum\limits_{l=0}^{k-1}\dfrac{\delta^{2l}}{l!^2}e^{\delta^2}
\ds\sum\limits_{\nu=0}^{\infty}\dfrac{1}{\nu!}\|f^{(l+\nu)}\|_2^2\nonumber\\
&\leq& \ds\sum\limits_{l=0}^{k-1}\dfrac{\delta^{2l}}{l!^2}e^{\delta^2}
\ds\sum\limits_{\nu=0}^{\infty}\dfrac{1}{\nu!}\sigma^{2(l+\nu)}\|f\|_2^2,\nonumber
\end{eqnarray}
using Bernstein's inequality. Thus
\begin{eqnarray}\label{pap3eqn2.14}
\ds\sum\limits_{i\in\mathbb{Z}}\ds\sum\limits_{l=0}^{k-1}|f^{(l)}(x_i)|^2c_{i,l}
&\leq& \left(\ds\sum\limits_{l=0}^{k-1}\dfrac{(\delta\sigma)^{2l}}{l!^2}\right)e^{\delta^2+\sigma^2}
\|f\|_2^2.
\end{eqnarray}
Similarly, we can prove that
\begin{eqnarray}\label{pap3eqn2.15}
\ds\sum\limits_{i\in\mathbb{Z}}\ds\sum\limits_{l=0}^{k-1}|f^{(l)}(x_i)|^2c_{i-1,l}
&\leq& \left(\ds\sum\limits_{l=0}^{k-1}\dfrac{(\delta\sigma)^{2l}}{l!^2}\right)e^{\delta^2+\sigma^2}
\|f\|_2^2.
\end{eqnarray}
From \eqref{pap3eqn2.14} and \eqref{pap3eqn2.15}, we get RHS of \eqref{pap3eqn2.7}.
\end{pf}

This leads us to the following frame algorithm: Recall that $B_\sigma$ is a reproducing kernel Hilbert space with reproducing kernel $K(x,y)=\dfrac{\sin \sigma\left(x-y\right)}{\sigma \left(x-y\right)}$.  \textit{i.e.,} every $f\in B_\sigma$ can be written as
\begin{eqnarray}
f(x)=\int\limits_{\mathbb{R}}f(t)\dfrac{\sin \sigma\left(t-x\right)}{\sigma \left(t-x\right)}~dt
=\langle f, K_x \rangle,
\end{eqnarray}
where $K_x(t)=K(t,x)$. Moreover, $f^{(r)}(x)=(-1)^r\langle f, K_x^{(r)} \rangle$. Hence, if $\delta<\dfrac{1}{\sigma}c_{k}^{\frac{1}{2k}}$,
it follows from Theorem \ref{pap3thm2.3} that the family $\{\sqrt{c_{i,l}+c_{i-1,l}}K_{x_i}^{(l)}:l=0,1,\dots, k-1, i\in\mathbb{Z}\}$ is a frame with frame bounds $A$ and $B$.\\\\
\noindent{}
\textsf{\underline{Frame Algorithm:}}

Set $S_kf:=\ds\sum\limits_{i\in\mathbb{Z}}
\ds\sum\limits_{l=0}^{k-1}(-1)^lf^{(l)}(x_i)(c_{i,l}+c_{i-1,l})K_{x_i}^{(l)}$
and $\rho=\dfrac{2}{A+B}$. Define
\begin{eqnarray*}
f_0&=&0,\\
f_{n+1}&=&f_n+\rho S_k(f-f_n),~ n\geq 0.
\end{eqnarray*}
Then we have $\lim\limits_{n\to\infty}f_n=f$. The error estimate after $n$ iterations turns out to be
\begin{eqnarray*}
\|f-f_n\|_2&\leq&\left(\dfrac{B-A}{B+A}\right)^n\|f\|_2.
\end{eqnarray*}
We refer to \cite{Raza} and also \cite{Grochenig1} for further details.

Now we shall provide a numerical comparison between the maximum gap conditions of \cite{Raza} and the current paper for various values of $k$.
Towards this end, we explicitly mention the values and bounds for the constants $c_r$ as given in \cite{BoWi}.
\begin{eqnarray*}
c_1=\pi^2, c_2=500.5467, c_3=61529.
\end{eqnarray*}
For any $r\geq 1$,
\begin{eqnarray}\label{pap3eqn2.18}
\dfrac{4r-2}{4r^2-r}\dfrac{(4r)!(r!)^2}{(2r!)^2}\leq c_r\leq\dfrac{4r+1}{2r+1}\dfrac{(4r)!(r!)^2}{(2r!)^2}.
\end{eqnarray}
Further, as $r\to \infty$, $c_r=\sqrt{8\pi r}\left(\dfrac{4r}{e}\right)^{2r}\left[1+O\Big(\dfrac{1}{\sqrt{r}}\Big)\right]$.
The upper bound  given in \eqref{pap3eqn2.18} is asymptotically exact. We take $\sigma=\pi$. 
\begin{center}
 \begin{tabular}{|l|l|l|}
  \hline
  \multicolumn{3}{|c|}{Maximum gap condition with $\delta<L$}\\
  \hline
  $k$ & ~~$L$ from \cite{Raza}~~~& $L$ from current paper\\
  \hline
  1 & 0.9003 &1\\
  \hline
  2 & 1.1849 & $\tau^{1/4}=1.5$ \\
  \hline
  3 & 1.4139 & 2\\
  \hline
  4 & 1.6479 & 1.9169 \\
  \hline
  5 & 1.8852&  2.3610\\
  \hline
  6 & 2.1239 & 2.8094\\
  \hline
  7 &2.3632 &3.2608\\
  \hline
  8 &2.6028 &3.7144\\
  \hline
  9 &2.8425 &4.1697\\
  \hline
  10 &3.0489 &4.6263\\
  \hline
  20 &5.4697 &10.0440\\
  \hline
  30 & 7.8448 &13.8689\\
  \hline
  \end{tabular}\\
\vspace*{0.3cm}
Table 1.1
\end{center}
In the Table 1.1, $L=\dfrac{2}{\pi}((k-1)!\sqrt{(2k-1)2k})^{1/k}$ from Razafinjatovo's method  and
$L=\dfrac{1}{\pi}c_{k}^{\frac{1}{2k}}$ from the current paper. We have made use of the lower bounds of $c_k$ given in \eqref{pap3eqn2.18} for $k\geq 4$ in the Table 1.1. Further if we  make use of the upper bound for $c_r$ given in \eqref{pap3eqn2.18}, we have the following table.
\begin{center}
 \begin{tabular}{|c|c|}
  \hline
  $k$& $L\approx$ \\
  \hline
  40 & 19.5623 \\
  \hline
  41 & 20.0333 \\
  \hline
  42 & 20.5043\\
  \hline
  \end{tabular}\\
\vspace*{0.3cm}
Table 1.2
\end{center}
Since the upper bound appearing in \eqref{pap3eqn2.18} is asymptotically exact, we expect that for sufficiently large $k$, the maximum gap may approach to $\dfrac{k-1}{2}$.

%

\end{document}